\documentstyle[12pt,seceqn]{article}
\font\BBB=msym10 scaled \magstep 1
\newtheorem{Lemma}{Lemma}[section]

\title { Julia sets and complex singularities \\
	in hierarchical Ising models}

\author{{\hbox{\vbox{\hsize=4.4in \baselineskip=10pt
	\hbox{P. M. Bleher \hfil} \vspace{5pt}        
	\hbox{\small School of Mathematical Studies \hfil}
	\hbox{\small Tel-Aviv University\hfil}
	\hbox{\small 69978 Israel\hfil}}}}
    {\hbox{\vbox{\hsize=4.4in \baselineskip=10pt
	\hbox{M. Yu. Lyubich\hfil} \vspace{5pt}	 
	\hbox{\small Institute for Mathematical Sciences\hfil}
	\hbox{\small SUNY, Stony Brook, \hfil}
	\hbox{\small N.Y. 11794, USA\hfil}}}}}
\date{March 4,1990}
\def\QED{\hfill{{\vrule height .7ex width .5em}}\par\medskip}

\input psfig
\begin{document}
\maketitle

\begin{abstract}
     We study the analytical continuation in the complex plane of free energy
of the Ising model on diamond-like hierarchical lattices. It is known 
\cite{DDI,DDL} that the singularities of free energy of this model lie on the
Julia set of some rational endomorphism $f$ related to the action of the Migdal-
Kadanoff renorm-group. We study the asymptotics of free energy when
temperature goes along hyperbolic geodesics to the boundary of an attractive
basin of $f$. We prove that 
for almost all  (with respect to the harmonic measure) geodesics
the complex critical exponent is common, and compute it.
\end{abstract}

\section {Introduction}

The purpose of this article is to analyse
complex singularities in temperature
of the free energy ${\cal F}$ in the Ising model
on diamond--like hierarchical lattices.
According to the traditional point of view a phase
transition manifests itself as a singularity of ${\cal F}$
as a function of thermodynamic parameters (like temperature
and external magnetic field).
 From this
point of view the theory of phase transitions should describe the
domain of analyticity of ${\cal F}$ and the type
of its singularities at points
of phase transition (see [1], where diverse approaches
to the first of these problems are discussed).

Since ${\cal F}$ is real analytic outside
of points of phase transition, it
can be continued into complex space with respect to
the thermodynamic parameters. Description of its complex singularities
is of great interest for the
theory of phase transitions because it
determines  analytic properties of the thermodynamic
function .

The celebrated Lee--Yang theory (see [2]) gives a
realisation of this approach describing the
 singularities of the analytic continuation
of the free energy in the ferromagnetic Ising model
with respect to the external magnetic field.
It proves that the zeroes of the grand partition function
in the ferromagnetic Ising model lie on the imaginary axis,
and hence complex singularities of the free energy
lie on the imaginary axis as well.
 An important problem stated in [2] is to study the limit distribution
of zeros of the grand partition function, since the free energy can be 
expressed as a logarithmic potential over this distribution.

The problem of description of complex singularities of the analytic
continuation of thermodynamic functions in temperature
is also very interesting from different points of view.
Many properties of asymptotic behavior of thermodynamic
functions in vicinity of a critical point were investigated
through the Kadanoff--Wilson--Fisher renormalisation group theory
(see e.g. [3,4]). It gives a local form of real critical singularities of
thermodynamic functions which has a nice
universal scaling structure. A problem is how are these singularities
continued to complex space and what is their global structure
in complex space?

 Unfortunately no general theory like the
Lee--Yang theory exists which describes for general models
global complex singularities
of thermodynamic functions in complex temperature plane .
However some exact results were obtained for the two dimensional Ising model.
The main tool here is the famous Onsager solution.
It turns out that in isotropic two dimensional Ising model
the zeroes of the partition function lie asymptotically
on two circles $e^{-2J/T}=\pm 1+\sqrt 2 e^{i\varphi}$ (it was conjectured
by Fisher  [5] and proved in [6]).
  Later it was shown (see [7,8] and references there)
that in anisotropic two dimensional Ising models on diverse
lattices the zeroes of the partition function fill some planar
regions in the complex temperature plane and in some cases
the density of the limit
distribution of zeroes can be found explicitly.

In  the present paper we consider another exactly solvable
model, namely the Ising model on diamond-like hierarchical lattices
(see [9-11]) . In this case the Migdal--Kadanoff
renormalization transformation turns out to be  to a rational map $f$ on the 
Riemann sphere which can be found explicitly.
The following nice observation was made in papers [12, 13]:
the set of complex points of phase transition coincides with the Julia set of
this map. In this paper we study the analytical properties of the
free energy {\cal F} near the singular points. We believe that clear
understanding of analytical properties of thermodynamical functions in the
complex plane has much to do with physical nature of the model.

The hierarchical sequence of the  diamond-like lattices
 depends on one natural parameter $b\geq 2$.
The lattice $\Gamma_0$ is just two "outer" sites  related by a bond.  In order
to obtain $\Gamma_1$ we insert in between the outer sites $b$ inner sites
related by bonds with outer ones (see Fig. 0 for $b=3$). Then in order to obtain 
$\Gamma_{n+1}$ we replace each bond in $\Gamma_n$ by the lattice $\Gamma_1$,
$n=1,2,...$.

\psfig {figure=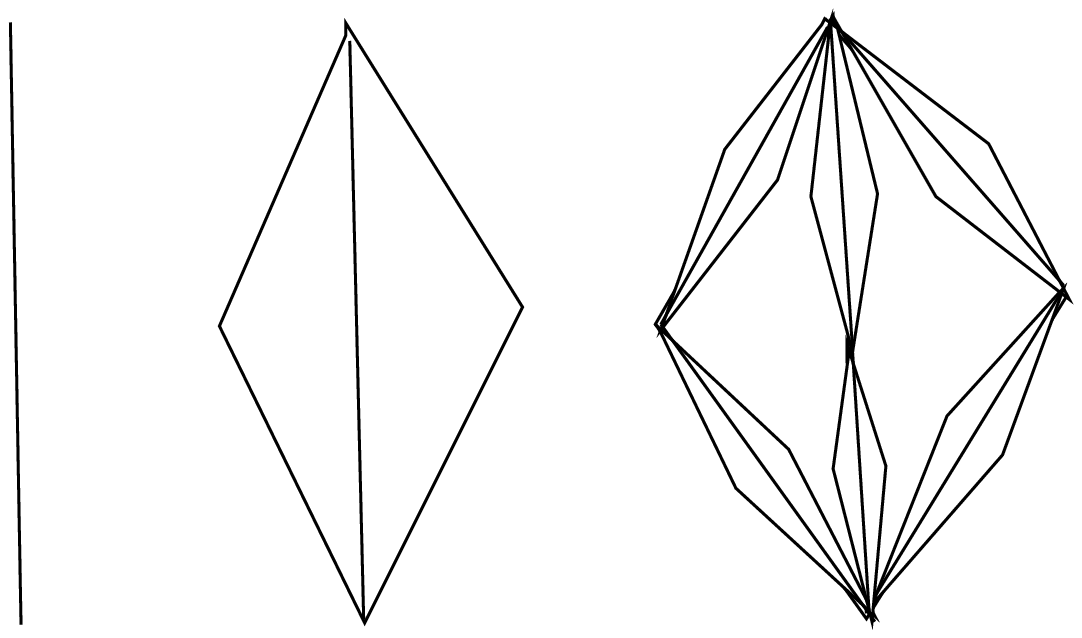,height=3in,width=4in}
\hspace{2in} Figure 0
\vspace{.5in}

We will refer to papers [11-15]  for the description of the Ising model
on these lattices and the calculation of the thermodynamical
functions. The starting point for us is the following
explicit formula for the free energy:

\begin{equation}
{\cal F}= -\frac{J}{2}-\frac{T}{2}\sum_{n=0}^{\infty}(2b)^{-n}\ln(1+t_n^b)
\label{energy}
\end{equation}
where $b=2^{d-1}$, $d$ is the ``dimension" of the lattice, $J$ is the interaction
constant, $T$ is temperature, and the sequence $t_n,\; n\geq 0$, is given by the
following recurrent equation

\begin{equation}
t_{n+1}=f(t_n),\qquad n\geq 0
\label{rec}
\end{equation}
where
\begin{equation}
f(t)=\frac{4t^b}{(1+t^b)^2},
\end{equation}
with the initial data
\begin{equation}
t_0=\exp(-\frac{2J}{bT})\equiv{G(T)}.
\end{equation} 

Equations (1.2)-(1.4) mean from the physical point of view that the map
$T\mapsto G^{-1}\circ f \circ G(T)$ gives the rescaling of temperature
under the Migdal-Kadanoff renorm-group transformation (see [9-15]). Note
that the points $t=0$ and $t=1$ are superstable fixed points of the map 
$t\mapsto f(t)$
(low- and high-temperature fixed points of the renorm-group) and for $b>1$ 
there exists the unique unstable fixed point $t_c$ on [0,1] . 
The critical temperature $T_c$ is equal to $G^{-1}(t_c)$.

     Formulas (1.1)-(1.4) make sense for complex values of $T$ as well. So,
we can consider the analytical continuation of free energy ${\cal F}(T)$ from
the positive axis $T>0$ into the complex plane. It is not hard to see
that the singularities of $\cal F$ lie on the Julia set
$J(f)$ [12, 13] 
(For the definition of the latter, see one of the surveys [7-10]).
     
     Let us consider now the immediate attracting basins $\Omega_0$ and 
$\Omega_1$ of points 0 and 1. One can show that $\overline{\bf C}\setminus{J(f)}$ is the union  of preimages of these domains, and $\Omega_0$ is a Jordan domain in $J(f)$ 
(see \S 2).In this paper we study the boundary properties of $\cal F$ in the
domain $\Omega_0$. To this end let us consider the Riemann map 
$\psi :\Omega_0\longrightarrow{\bf U}$ of $\Omega_0$ onto the unit disk. 
The hyperbolic geodesics in $\Omega_0$ are just the $\psi^{-1}$-
images of the radii in ${\bf U}$. 
Denote by $B_\tau$ the geodesic ending at $\tau\in\partial\Omega_0$.  
Let us consider also the harmonic measure $\mu$ on $\partial\Omega_0$, i.e.
${\mu}={\psi_{\ast}^{-1}\lambda}$, where $\lambda$ is the Lebesgue measure on the
circle ${\partial}{\bf U}{\equiv}{\bf T}$. For $t\in{B_{\tau}}$ denote by 
$l(t)$ the length
of $B_{\tau}$ from $t$ to $\tau$ (perhaps, $l(t)=\infty$).
    
      In the present paper we prove:

(i) The derivative ${\cal F}^{\prime}$ of free energy is continuous up to the
boundary of $\Omega_0$.

(ii) For $b > 2$ the second derivative is discontinuous in $\Omega$, and
has the following asymptotics on $\mu$-almost all geodesics:

\begin{equation}
 \lim_{t\rightarrow\tau , t \in{B_\tau}}  \frac{\ln \mid {\cal F}''(t)\mid}{-\ln l(t)}\equiv\alpha_c=1-\frac{\ln 2}{\ln b}=1-\frac{1}{d-1}.
\end{equation}

This means that for almost all geodesics the specific 
heat critical exponent
in the region of low temperatures is universal and equal to
$1-\frac{1}{d-1}$. 

     Now let us dwell in more detail on the content of the paper. 

In  \S 2 we
describe the dynamical properties of the endomorphism $f$. In particular, we
show using the Douady-Hubbard straightening theorem that $\Omega_0$ is a
Jordan domain.

    In   \S 3 we show that ${\cal F}'$ is continuous in $\rm cl\Omega_0$
and that $\partial\Omega_0$ is the natural boundary of analyticity of ${\cal F}$.
The proof is based upon some amusing observations concerning $f$ (its relation
to the Koebe function and a Tchebyshev polynomial).

    In  \S 4 we discuss some
technical background: the Bowen-Ruelle-Sinai  thermodynamical formalism and
the construction of the natural extension (the inverse limit) of $f$. These
are the main tools (together with the ergodic theorem) for the accurate 
computation of the critical exponent. 

     In \S 5 we discuss the functional equation for ${\cal F}''$ and related
spectral properties of the weighted substitution operator in the disk-algebra.

    \S 6 is the central section of the paper: here we
give the computation of the critical exponent, provided ${\cal F}''$ is not
continuous up to the boundary of $\Omega_0$. 

In  \S 7 we prove that
${\cal F}''$ really satisfies this property, which completes the proof of
the main result. 

     In the last \S 8 we discuss some related problems.

{\bf Acknowledgement.} We are grateful to J. Milnor for looking through the text 
and making useful remarks,  to M. Fisher and the referee for critical
comments yielding the improvement of the exposition,
 to J. Milnor and G. Tusnagy for  making nice computer pictures, 
and to NEFIM fund of
Hungarian Academy of Sciences for the support of the visit of one of the authors
(P.M.B.) to Budapest. 
     
\section{Dynamics of the map $f:t\mapsto\frac{4t^b}{(1+t^b)^2}$}

     We refer to the surveys [16-19] for the general view of the dynamics of
complex rational maps. We will use some concepts and facts of this theory
without extra explanations.

     Let us introduce the following notations:\\
$f^{\circ n}=f\circ...\circ f$ is the $n$-fold iterate of $f$;\\
$C(f)$ is the set of its critical points (a rational map of degree $d$ has
     $2d-2$ \hspace{.5in} critical points counting with multiplicity);\\
$\mid\mid\cdot\mid\mid$ is the spherical metric on $\overline{\bf C}$;\\
${\bf U}=\{z:\mid z \mid\leq 1\} $ is the closed disk;\\
${\bf U}^\circ={\rm int}{\bf U}$ is its interior;\\
${\bf T}=\partial {\bf U}$ is the unit circle;\\
$B(a,r)=\{z\in {\bf C}:\mid z-a\mid\leq R\}$ for $a\in {\bf C}$;\\
$J(f)$ is the Julia set of $f$. 

    The function 
\begin{displaymath}
   f\equiv f_b :t\mapsto \frac{4t^b}{(1+t^b)^2}
\end{displaymath} 
is related to the well-known {\em extremal Koebe function}(see [20])
\begin{displaymath}
    {\cal K}_0 (z)=\frac{z}{(1-z)^2}.
\end{displaymath}
Setting $ K(z)=-4{\cal K}_0 (-z) $ and $S(z)\equiv S_b (z)=z^b$
we have $f=K\circ S$.

     The relation of $f$ to the Koebe function is quite mysterious, especially
if one relates the coefficient 4 to the Koebe constant 1/4. It becomes still
more amusing if to observe that $K(t)$ is conformally conjugated to the
{\em Tchebyshev polynomial} $T:\tau\mapsto 2{\tau}^2-1$. Indeed, the function $K$
has two simple critical points $c_1=1$ and $c_2=-1$. Moreover, $f(c_1)=c_1$,
i.e.\ $c_1$ is superstable fixed point, and $c_2\mapsto\infty\mapsto 0$,
where 0 is repelling fixed point. Up to conformal conjugation, $T$ is the 
unique rational function of degree 2 possessing such properties. More 
specifically, $\varphi\circ K\circ {\varphi}^{-1}=T$ where 
$\varphi :t\mapsto\frac{1+t}{1-t}$ is the M\"{o}bius transformation mapping
the triple $\{0,1,\infty\}$ onto the triple $\{1,\infty ,-1\}$. In particular,
it follows that the Julia set $J(K)$ coincides with the negative semi-axis
$[-\infty ,0]={\varphi}^{-1} [-1,1]$.
 
     The power functions and Tchebyshev polynomials play a particular role
in the iteration theory. They appear as the exceptions in a number of problems;
e.g.,\ only these functions have a Julia set with simple geometry. The 
composition $f=K\circ S$ does not possess such a property(see Fig. 1).
      
        A rational function $g$ is called {\em critically finite} if the orbits
\{$g^n(c_i)$\} of all its critical points are finite.

       By the chain rule
\[ C(f)=C(S)\cup S^{-1} C(K)=\{0,\infty\}\cup\{\alpha_i\}_{i=1}^b\cup\\             \{\beta_j\}_{j=1}^b ,\]
where $\alpha_i$ are the $b$th roots of 1 , and $\beta_j$ are 
the $b$th roots of -1.
Moreover, 0 and $\alpha_1=1$ are superstable fixed points which absorb the      
 orbits of all other critical points: $\alpha_i\mapsto 1$,
$\beta_j\mapsto\infty\mapsto 0$. Thus, the function $f$ is critically finite.

      Denote by $\Omega_a\equiv\Omega(a)$ , the component
of $F(f)$ containing $a$. The domains
$\Omega_0$ and $\Omega_1$ are called the {\em immediate basins} of the fixed
 points 0 and 1.

      We say that a rational function $g$ {\em satisfies the axiom $A$}
if the following equivalent properties hold :

(i) The orbits of all critical points converge to stable cycles;

(ii) $g$ is expanding on the Julia set,i.e.,\ there exist constants $C>0$
     and $\lambda>1$ such that
\[\mid\mid dg^n(z)\mid\mid\geq C\lambda^n (z\in J(f), n\in {\bf N}).\]
       
       It follows from above that our function $f$ satisfies (i) and hence
satisfies axiom $A$. This implies in particular that the Fatou set
consists of the preimages of the immediate basins $\Omega_0$ and $\Omega_1$.

     Set $\Omega=\rm cl\Omega_0$ and $\Gamma=\partial\Omega_0$ 
(these notations will  be used up to the end of the paper).
 We will show now that $\Gamma$ is a
Jordan curve and even a quasicircle\footnote{actually, the last holds  
          automatically} (but by Fatou's theorem it has no tangents
at any point). To this end we apply the Douady-Hubbard  straightening theorem
(see [21]).   

       Let $V$ and $V'$ be two simply connected domains bounded by         
      piecewise-smooth curves, and ${\rm cl}V\subset V'\subset{\bf C}$. A map
$g:V\rightarrow V'$ is called {\em polynomial-like of degree $d$} if it is a
$d$-sheeted analytical covering of $V$ over $V'$ having no critical points on
$\partial V$. By the Riemann-Hurwitz formula, such a map has $d-1$ critical     points in $V$ counting with multiplicity.  Set

\[ K(g)=\{z:g^n(z)\in V (n=0,1,...)\}, K^0(g)={\rm int} K(g).\]
The $K(g)$ is a compact subset of $V$.

           The straightening theorem states that any polynomial-like map $g$
is quasi-conformally conjugated to a polynomial of the same degree, i.e.,\
there exists a quasi-conformal homeomorphism $\psi :{\bf C}\rightarrow{\bf C}$
such that $\psi\circ g\mid W=h\circ\psi\mid W$ for some domain $W$, 
$K\subset W\subset V$.
 Moreover, $\psi\mid K^0 (g)$ is conformal and
  \[\psi (K(g))={\bf C}\verb+\+ \{z: h^{\circ n}(z)\rightarrow\infty\                             (n\rightarrow\infty )\}\]
is the {\em filled-in Julia set of $h$}.

\newtheorem{Jordan}{Lemma}[section]. 

\begin{Jordan}
The domain $\Omega_0$ is Jordan, and its boundary is a
quasi-circle. The restriction $f\mid\Omega_0$ is conformally conjugated
to the power transformation $z\mapsto z^b$ of the unit disk ${\bf U}^0$.

\end{Jordan}
{\bf Proof}. Let us construct a neighbourhood of $\Omega$ on which $f$ is a
polynomial-like map of degree $b$. To this end note that the function $K$
conformally maps the disk ${\bf U}^0$ onto the plane slitted along
the semi-axis ${\bf R}_1=[1,\infty ) $                                             (this is the characteristic property of the 4-fold
Koebe function). Let us consider the domain $V$(see Fig. 2) bounded
by the arc ${\gamma}_1$ of the circle $B(1,\varepsilon )$, the arc ${\gamma}_2$
of the circle $B(0,R)$ where $R>>1$ and two horizontal intervals. 

\psfig {figure=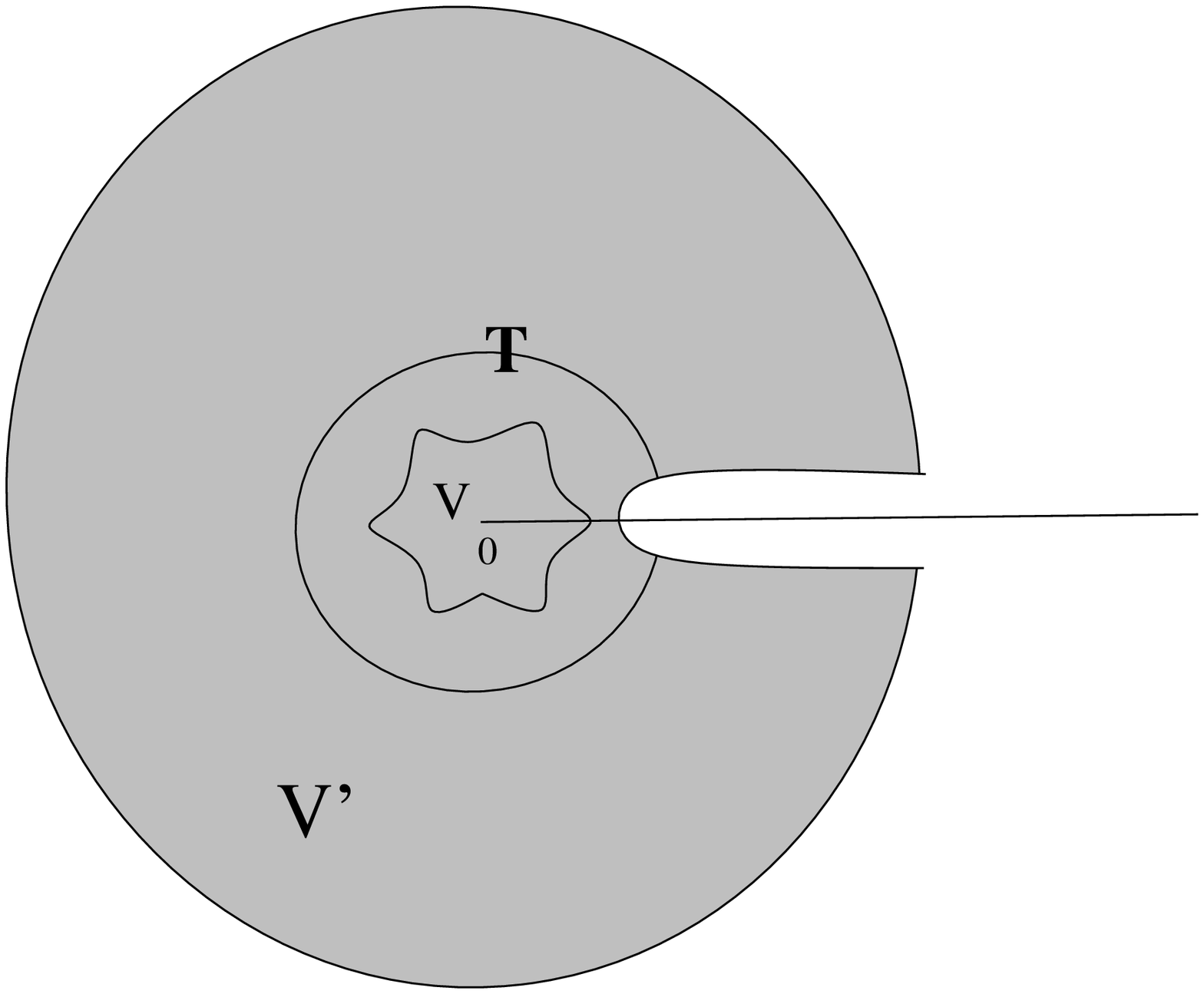,height=3in,width=4in}
\hspace{2in} Figure 2
\vspace{.5in}

      Let $V$ be the component of the inverse image ${f^{-1}}(V') $
containing 0. Then ${\rm cl}V\subset {\bf U}^0 $ since 
$f(\partial{\bf U})={\bf R}_1$ lies outside $V'$. Besides, 
${\rm cl}V\cap{\gamma}_1 =\emptyset$ for sufficiently small $\varepsilon$. Indeed, 
as 1 is a stable fixed point, the arc $f({\gamma}_1)$ lies inside the disk
$B(1,\varepsilon)$ and, hence, outside $V'$.

       Thus, ${\rm cl}V\subset V'$.

     Further, it is clear from $V=(K\mid{\bf U})^{-1} \circ (S^{-1} V')$
that $V$ is simply-connected. Indeed, it is elementary that $S^{-1} V'$ is
simply-connected (see e.g.,\ [17], Lemma 1.4), 
while $K\mid {\bf U}$ is univalent.

     We have shown that $f:V\rightarrow V'$ is a polynomial-like map.
Its degree is equal to $b$ since $V\subset {\BBB U}_0$ contains the unique
$(b-1)$-fold critical point 0. Clearly, $\Omega\subset K(f\mid V)$.

        By the Straightening Theorem, $f:V\rightarrow V'$ is quasi-conformally
conjugated to a polynomial $h$ of degree $b$. Normalize the conjugating 
homeomorphism $\psi$ in such a way that $\psi (0)=0$. Then 0 becomes a
$(b-1)$-fold critical point for $h$. It follows that $h(z)=Cz^b$. Normalizing
$\psi$ additionally in such a way that $\psi (t_c)=1$ where $t_c$ is a real
fixed point lying on $\Gamma$ ,we get $C=1$. 

      Thus, $\psi$ conjugates $f: V\rightarrow V'$ to the power polynomial
$h: z\mapsto z^b$. Consequently, $\Omega_0 =\psi^{-1}({\bf U})$ is a Jordan
domain bounded by a quasi-circle, and $\psi$ conformally conjugates 
$f\mid\Omega_0$ to $h\mid{\bf U}$. The lemma is proved.  \QED

{\bf Remark.} Tan Lei showed us another proof of the above lemma which can be
applied to the high-temperature region as well. 

\section{Analytic properties of ${\cal F}'$}

     In this section we will show that the derivative ${\cal F}'(t)$ is continuous
in the closed set $\Omega$. So, ${\cal F}'$ belongs to the disk-algebra 
$A(\Omega)$,
i.e.\ the  algebra of functions continuous  in $\Omega$ and holomorphic in
$\Omega_0$. We will get this estimating the derivative $\mid\mid Df\mid\mid_\nu$
in a special Riemann metric $\mu$.

      From now on we will consider the following function $F$ instead of free
energy ${\cal F}$ (1.1):

\begin{equation}
 F=\sum_{n=0}^\infty \frac{1}{(2b)^n}g\circ f^n
\end{equation}
where $g(t)=\ln (1+t^b)$. Clearly, its analytical properties are the
same as those of $\cal F$. Note that $g$ is analytic in a neighbourhood of
$\Omega$, and so series (3.1) converges uniformly in $\Omega$. Hence
$F\in A(\Omega)$. 
       Further, for $z\in\Omega_0$ we have

\begin{equation}
 F'(z)=\sum_{n=0}^\infty \frac{1}{(2b)^n} g'(f^n(z)) (f^n)'(z)
\end{equation}

       We want to show that this series converges uniformly in $\Omega$,
which certainly implies $F'\in A(\Omega)$. The required statement follows
from the following estimate:

\begin{Lemma}
   $\mid (f^n )'(z)\mid\leq C(\sqrt{2} b)^n$ for $z\in\Omega$.
\end{Lemma}        
   {\bf Proof}   
Let us recall that $f=K\circ S$  (see \S 2). The power function $S$ satisfies
the functional equation $\exp (bz)=S(\exp z)$. From the dynamical viewpoint
it means that $\exp$ semiconjugates the transformations 
$L:z\mapsto bz$ and $S$. Denote by $\sigma $ the Euclidean  metric on 
${\bf C}$, and by $\mu =\exp_*\sigma$ its image on the punctured plane
${\bf C}^* ={\bf C}\verb+\+\{0\}$. We have 
 $\mid d\mu\mid=\mid dt\mid/\mid t\mid$.

      As $\mid\mid DL(z)\mid\mid_\sigma =b$ for $z\in {\bf C}$,

\begin{equation}
     \mid\mid DS(t) \mid\mid_\mu =b    
\end{equation} 
     for $t\in{\bf C}^*$. Besides

\begin{equation}
\mid\mid DK(t) \mid\mid_\mu =\frac{\mid tK'(t) \mid}{\mid K(t) \mid}=
   \frac{\mid 1-t\mid}{\mid 1+t \mid}=\frac{1}{\mid \varphi (t) \mid}
\end{equation}
  where $\varphi (t)=(1+t)/(1-t)$. It is surprising that exactly this function
conjugates $K$ to the Tchebyshev polynomial $T:\tau\mapsto 2\tau^2 -1$.
Due to this observation it is reasonable to pass to the conjugated function
    \[h=\varphi\circ f\circ\varphi^{-1} =T\circ R\]
where $R=\phi\circ S\circ \phi^{-1}$. Consider the corresponding Riemannian
metric $\nu =\phi_*\mu$. By (3.3), (3.4),
      \[\mid\mid DR(\tau)\mid\mid_\nu =b, \quad
\mid\mid DT(\tau)\mid\mid_\nu =\frac{1}{\mid\tau\mid}.\]
      Hence for $\tau =\varphi (t)$
\begin{equation}
     \mid\mid Df(t) \mid\mid_\mu =\mid\mid Dh(\tau) \mid\mid_\nu =
     \frac{b}{\mid R(\tau)\mid}.    
\end{equation}

    The function $h$ has a superstable fixed point $1=\varphi(0)$ with the
immediate attracting basin $W^0=\varphi\Omega_0$. Since 
$\Omega\subset{\bf U}^0$, the set $W\equiv {\rm cl}W^0 =\varphi(\Omega)$ lies in
the right half-plane $\varphi({\bf U}^0)={\bf P}^0=\{\tau:Re\tau>0\} $.
As $W$ is $g$-invariant, $T(R(W))\subset {\bf P}^0$ and hence (see Fig. 3)

\[ R(W)\subset T^{-1}({\bf P}^0) =\{\tau=x+iy: x^2-y^2>1/\sqrt{2}\}\]

Consequently, $\mid R(\tau)\mid>1/\sqrt{2}$ for $\tau\in W$, and (3.5) implies
    
       \[\mid\mid Df(t) \mid\mid_\mu < b\sqrt{2},\quad t\in\Omega \]

Hence,  \[\mid\mid Df^n (t) \mid\mid_\mu < (b\sqrt{2})^n,\quad n\in{\bf N} .\]
But on the boundary $\partial\Omega$ the metric $\mu$ is equivalent to the
Euclidean metric and, hence, $\mid (f^n)'(t)\mid <C(\sqrt{2} b)^n$
for $t\in\partial \Omega$. By the Maximum Principle, this inequality holds for
$t\in\Omega$. The lemma is proved. \QED

\psfig{figure=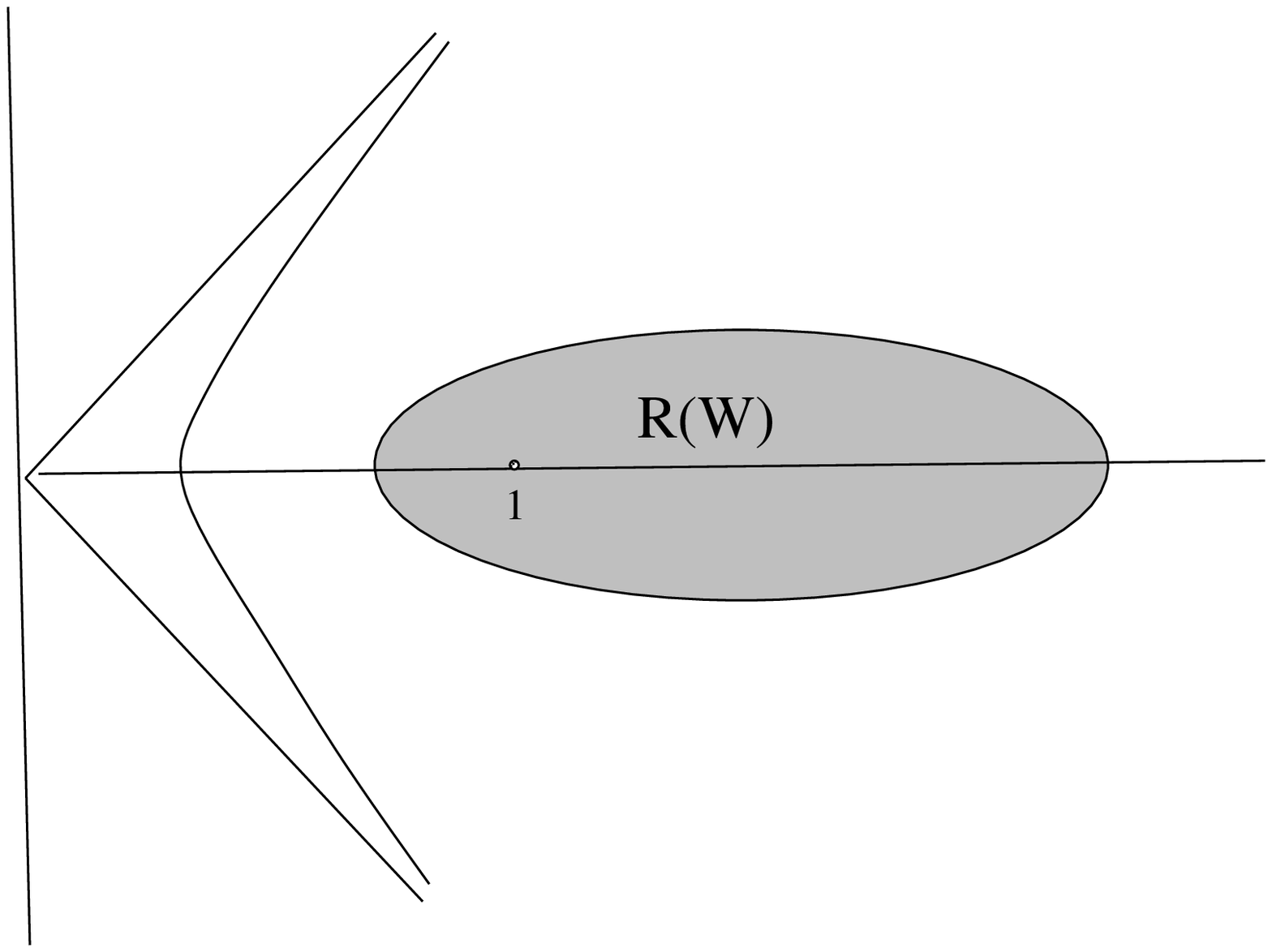,height=2.5in,width=3in}
\hspace{2in}  Figure 3
\vspace{.5in}

       Now let us establish some global analytical properties of $F'$. It is
curious that this function has no singularities at points 
$\beta_i=\sqrt[b]{-1}$.

\begin{Lemma}
{\rm (i)} The function $F'$ is a single-valued holomorphic function
on the Fatou set $N(f)$.
      
{\rm (ii)}The set $\Omega^0$ is the maximal domain of analyticity of $F'$.
\end{Lemma}
{\bf Proof}. Consider the multy-valued function
$\sigma(t)=g(t)+\frac{1}{2b} g(f(t))$ where $g(t)=\ln(1+t^b)$ as above.
Let us show that it is regular near  $\beta_j =\sqrt[b]{-1},\quad j=1,...,b$.
We have
  \[ \sigma(t)=\ln(1+t^b)+\frac{1}{2b}\ln(1+\frac{4^b t^{b^2}}{(1+t^b)^{2b}})=\]

         \[  =\frac{1}{2b}\ln((1+t^b)^{2b}+4^b t^{b^2}) .\]

So, $\sigma(\beta_j)=\frac{1}{2}\ln4+(\frac{\pi}{2}+ \pi n)i$, and we see
that $\beta_j$ are regular points for $\sigma$.  Hence,$\quad\sigma$ is regular
in the components $\Omega(\beta_j)$.        
       
         Now let $V$ be an $n$-fold preimage of some $\Omega(\beta_j)$. By (3.1),

    \[ F(t)=\frac{1}{(2b)^n}\sigma(f^n t)+ regular function.\]
Consequently, $F$ is regular in $V$.

       Thus, $\Omega(\infty)$ is the only component in which $F$ is not regular,
and  there we have 
    \[ F(t)=\ln(1+t^b)+ regular function.\]
Hence,
     \[ F'(t)=\frac{bt^{b-1}}{1+t^b}+ regular function\]
is  regular  in  $\Omega(\infty)$.

       So, $F'$ is regular on the whole Fatou set $N(f)$,and it is obvious
from (3.2) that it is single-valued. The (i) is proved.
       
        (ii) Let us consider the functional equation for $F$:
\begin{equation}
                  \frac{1}{2b}F\circ f-F=g
\end{equation}
 
   Taking the derivative, we get           
\begin{equation}
             \frac{f'}{2b}F'\circ f-F'=g'=\frac{bt^{b-1}}{1+t^b}.
\end{equation}
    Provided $F'$ can be analytically continued beyond $\Omega^0$ into some
neighbourhood $U$ of $t\in\partial\Omega$ , it follows from (3.7) that it
can be continued as a meromorphic function into $f^n U$. Since
$f^n U\supset J(f)$ for some $n$ , the function $F'$ is meromorphic on
 the whole sphere, i.e.\ rational.

       But we have shown that $F'$ has no poles in $N(f)$. If $F'$ had a pole
at a $t\in J(f)$ then by (3.7) it would have had poles at all points of the 
grand orbit
  \[ {\cal O}(t)= \bigcup_{m=0}^{\infty}\bigcup_{n=0}^{\infty} f^{-n}(f^m t)\]
But it is impossible since this orbit is infinite. So, $F'$ has no poles at all,
and this absurdity completes the proof. \QED

\section {Symbolic dynamics, thermodynamical formalism and the natural
  extension of $f\mid\Omega$}

   To study the boundary behaviour of $F''$ we need the Bowen-Ruelle-Sinai
thermodynamical formalism   (see[22],[24] or [18]). We will state
 the main results of
this theory for our particular map $f\mid\Gamma$. The theory can be applied
to this map without any problem because it is expanding (see \S 2).

      Let us consider the homeomorphism $\psi: \Omega\rightarrow{\bf U}$
conformal on $\Omega_0$ and reducing $f\mid\Omega$ to $S:z\mapsto z^b$
(Lemma 2.1). Let $\psi(t)=re^{2\pi i\theta}\in{\bf U}$. Associate to a point
$t\in\Omega$ the pair $(r,\overline\varepsilon_+)$ where $r\in[0,1]$ , and
$\overline\varepsilon_+ =(\varepsilon_0,\varepsilon_1,...)$
is the $b$-adic decomposition of the $\theta\in [0,1)$.  Then $f$  turns
into the transformation 
$(r,\overline\varepsilon_+)\mapsto(r^b,\sigma_+\overline\varepsilon_+)$
where 
\[\sigma_+:(\varepsilon_0,\varepsilon_1,...)\mapsto\\
           (\varepsilon_1,\varepsilon_2,...) \]
is the shift on the space $\Sigma_b^+$ of all one-sided $b$-adic sequences.
Sometimes we will identify $\Omega$ with $[0,1]\times\Sigma_b^+$ ,though the
described correspondence is not one-to-one. Then $f\mid\Gamma$ will be
identified with the shift $\sigma_+$.

     Denote by $[\varepsilon_0,...\varepsilon_{n-1}]$ $b$-adic cylinders
in $\Sigma_b^+$ and corresponding $b$-adic intervals in $\Gamma$. Let
$B_\theta =\psi^{-1}\{re^{2\pi i\theta}: 0\leq r\leq 1\}$ be {\em hyperbolic
geodesics} in $\Omega$, and 
$\Gamma_r =\psi^{-1}\{re^{2\pi i\theta}: 0\leq \theta\leq 1\}$ 
be {\em equipotential levels} (for the Green function $\ln\mid\psi\mid$).

     Let $\rho$ be a H\"{o}lder function on $\Gamma$ which is called the 
potential (here we pass from electrostatics to thermodynamics).Set
          \[S_n\rho =\sum_{k=0}^n{\rho\circ f^k}.\]
The Gibbs measure $\nu_\rho$ on $\Gamma$ corresponding to the potential
$\rho$ is the measure satisfying the following estimates on cylinders:

\begin{equation}
 \nu_\rho[\varepsilon_0...\varepsilon_{n-1}]\asymp
 \exp[S_n\rho(t_{\varepsilon_0...\varepsilon_{n-1}}) -nP],
\end{equation}          
where $t_{\varepsilon_0...\varepsilon_{n-1}} $ is any point of
$[\varepsilon_0...\varepsilon_{n-1}],$ and $P=P_f(\rho) $ is a constant
called {\em the pressure}, and the sign  $"\alpha\asymp\beta"$ means
$C_1\beta\leq\alpha\leq C_2\beta$.

      The main result of the Bowen-Ruelle-Sinai theory states that for any
H\"{o}lder function $\rho$ there exists the unique Gibbs measure $\nu_\rho$.
This measure satisfies {\em the Variational Principle}
\begin{equation}
  \sup_{\nu\in M(f)} (h_\nu(f) +\int_\Gamma{fd\nu})=
      h_{\nu_\rho}(f) +\int_\Gamma{fd\nu_\rho}= P_f(\rho)
\end{equation}
where $M(f)$ denotes the compactum of all $f$-invariant probability measures
on $\Gamma$, and $h_\nu (f)$ is the entropy of $\nu$.

     The pressure $P_f(\rho)$ is the smooth convex functional of $\rho$, and
its differential at $\rho$ is the Gibbs measure $\nu_\rho$ (see [23]):
\begin{equation}
  \frac{dP_f(\rho+\kappa\alpha)}{d\kappa}\mid_{\kappa =0} =
         \int{\alpha d\nu_\rho}
\end{equation}

       To the potential $\rho =0$ corresponds the unique measure of maximal
entropy $\mu\equiv\nu_0$, namely the Bernoully measure with $b$ equal states
(in view of the model $\Sigma_b^+$). The entropy of this measure is equal to the
topological entropy of $f\mid\Gamma$ : $h_\mu(f)=h(f)=\ln b$.

    The Riemann map $\psi :\Omega\rightarrow{\bf U}$ transforms $\mu$ into
the measure of maximal entropy for $z\mapsto z^b$, i.e.\ to the Lebesgue
measure on ${\bf T}$. Hence, $\mu$ coincides with the harmonic measure on
$\Gamma$ corresponding to 0 (see [20]). Consequently,
    \[H(0)=\int_\Gamma H d\nu \]
for any function $H$ harmonic in $\Omega_0$ and continuous up to the
boundary.

      Applying this formula to the function 
\footnote {taking in account that $\psi (t)\sim t$, since there is the following
formula for the Green function (cf.[32]):
 \[\ln\mid\psi\mid=\lim_{n\rightarrow\infty}d^{-n}{\mid f^{\circ n}\mid}\]}.
     \[H(t)=\ln\frac{\mid f'(t)\mid }{\mid\psi(t)\mid^b},\]
we find {\em the characteristic exponent} of $\mu$
\begin{equation}
       \chi_\mu =\int_\Gamma{\ln\mid f'\mid d\mu}=\ln b.
\end{equation}

      Let us pass now to the crucial construction of {\em the natural extension} 
(or {\em the inverse limit}) 
  $\overline f:\overline\Omega\rightarrow\overline\Omega$ (see[25]).
By definition, a point $\overline z\in\overline\Omega$ is the inverse
orbit  $\overline z =(z_0, z_{-1},...)$, i.e.\ $fz_{-(i+1)}=z_{-i} $, and
$\overline f :\overline z\mapsto(fz_0,z_0,z_1,...)$. The transformation
 $\overline f$ is invertible on $\overline\Omega$, and there exists the
natural projection $\pi:\overline\Omega\rightarrow\Omega$, $\pi(z)=z_0$,
semiconjugating $\overline f$ and $f$. All fibers of $\pi$ except zero one
are Cantor sets. Each $f$-invariant measure $\nu$ on $\Gamma$ can be uniquely
lifted to the $\overline f$-invariant measure $\overline\nu$ on 
$\overline \Gamma$.

     The symbolic dynamics for $f$ generates the symbolic dynamics for 
$\overline f$.
Namely, $\overline\Omega$ can be identified mod 0 with 
$[0,1]\times \Sigma_b$ where
 $\Sigma_b =\{(...\varepsilon_{-1},\varepsilon_{0},\varepsilon_{1}...)\}$
is the space of two-sided sequences. Then
    \[\overline f:(r,\overline\varepsilon)\mapsto    
     (r^b,\sigma\overline\varepsilon),\]
where $\sigma:\Sigma_b\rightarrow\Sigma_b$
is the left shift, and $\overline\mu$ turns into the Bernoulli measure on
$\Sigma_b$ with equal states.

       The lifts $\overline\Gamma_r =\pi^{-1}\Gamma_r$ of the equipotential   
levels will be called the {\em solenoids},
 $\overline\Gamma\equiv\overline\Gamma_1$. The reason is that these sets can
be supplied with the structure of the solenoidal group ${\bf T}$ ,i.e.\
the inverse limit of the group endomorphism $z\mapsto z^b$ of ${\bf T}$
 (see [26]).
Then $\overline f\mid\overline\Gamma$ turns into a group endomorphism,
and $\overline\mu$ into the Haar measure on $\overline\Gamma$.

     The space $\overline\Omega$ can be regarded as a continuum-sheeted
Riemann surface over $\Omega$, the "bunch of sheets" gluing together at
zero. If one cuts $\Omega$ along the geodesic $B_0,\quad  \overline\Omega$
is folliated into the sheets $L(\overline\varepsilon_-)$ coded by
one-sided sequences 
$\overline\varepsilon_- =(...\varepsilon_{-2},\varepsilon_{1})\in\Sigma_b^-$.

The gluing of the sheets is fulfilled  by the $b$-adic shift
$A\equiv A_b:\overline\varepsilon_-\mapsto\overline\varepsilon_- +{\bf 1}$
where {\bf 1}=(...0,0,1), and addition is understood in the sense of the
group of $b$-adic numbers . Gluing together countably many sheets
corresponding to an orbit $\{A^n(\overline\varepsilon_- )\}_{n=-\infty}^\infty$
of the $b$-adic shift, we get the logarithmic Riemann surface 
$W(\overline\varepsilon_-)$ \footnote{The inner topology of this surface
differs from the topology induced from the space $\overline\Omega$ in which
it is densely immersed}. All inverse functions $f^{-n}(z)$ become
single-valued on this surface.

     From the dynamical point of view, Riemann surfaces 
$W(\overline\varepsilon_-)$ are {\em global unstable manifolds} of $f$ :
if $\overline t,\overline\tau\in W(\overline\varepsilon_-)$, then 
$\overline f^{-1}$-orbits of   $\overline t$ and $\overline\tau$ are 
exponentially drawing together.

    For a function $\rho(t)$ in $\Omega$ we will write 
$\rho(\overline t)\equiv\rho(t)$, where $t=\pi(\overline t)$; in particular,
 $\mid\overline t\mid\equiv\mid t\mid$. Set
\begin{equation}
   W_\delta(\overline\varepsilon_-)=\{ \overline t\in W(\overline\varepsilon_-)      :\mid\overline t\mid\geq\delta\},\\
  S_{-n}\rho(\overline t)=-\sum_{k=0}^{-n+1}\rho(f^k(\overline t)).
\end{equation} 
The  drawing together of inverse orbits originating in $W_\delta(\overline\varepsilon_-)$
is uniformly exponential. It follows in a standard way that for any function
$\rho(t)$, H\"{o}lder on $\Omega^* =\Omega\verb+\+ \{0\}$, for $\delta>0$
and  $\overline t,\overline\tau\in W_\delta(\overline\varepsilon_-)$
the following estimates hold
\begin{equation}
   \mid S_{-n}\rho(\overline t)- S_{-n}\rho(\overline \tau)\mid\leq C(\delta),
\end{equation}
   where $C(\delta)$ does not depend on $n$ and $\overline\varepsilon_-$.

\section{ The weighted substitution operator}

      Let us consider a function $\beta\in A(\Omega)$ having no zeroes on
$\Gamma\equiv\partial\Omega$, and $h\in A(\Omega)$. In this and the next
sections we will consider the following functional equation:

\begin{equation}
     \beta(t)U(ft)-f(t)=-h(t)
\end{equation}

Let us consider the multiplicative cocycle 
$\beta_n(t)=\beta(t)\beta(ft)...\beta(f^{n-1}t)$
associated with $\beta$.    
 We  will assume in what follows that

(i)the cycle $\beta$ has the positive Liapunov exponent: 
\begin{equation}
   \chi_\mu(\beta)\equiv\int{\ln\mid\beta\mid} d\mu >0
\end{equation}

(ii)\footnote{this assumption is convenient but it is not essential}
The function \, $\ln\!\mid\beta\mid$ \, on $\Gamma$ is not 
homologous to a constant. 
 This means that there are no continuous solutions 
$\varphi\in C(\Gamma)$ of the equation
\[\ln\mid\beta\mid=\varphi\circ f-f+c\] with any constant $c$.
   
 {\bf Motivation.}  By differentiating equation (3.6), we get the equation 
(5.1) for $U=F''$ with:

\begin{equation}
  \beta(t)=\frac{1}{2b} f'(t)^2 ,\quad h(t)=\frac{1}{2b}F'(ft)f''(t)+g''(t).
\end{equation}
    The function $\beta$ is holomorphic in a neighbourhood of $\Omega$
and has there the unique root at $t=0$, and  the function $h\in A(\Omega)$
(since $F'\in A(\Omega)$). Assumption (i) is valid due to (4.4):

\begin{equation}
   \chi_\mu(\beta)=2\chi_\mu - \ln 2b=\ln\frac{b}{2} > 0. 
\end{equation}
Assumption (ii) holds since $\ln\mid f'\mid$ is not homologous to a constant 
on $\bf T$ [30].    
                                                             \QED

 The equation (5.1) has the unique solution holomorphic
in $\Omega_0$ :

\begin{equation}
    U(t)=\sum_{n=0}^\infty \beta_n (t) H(f^n t),
\end{equation}
where $\beta_n(t)$ is the multiplicative cocycle generated by the 
function $\beta$.

    In this section we will explain why this solution is not, as a rule,
continuous up to the boundary. The words ``as a rule" means: for any
$h\in A(\Omega)$ outside a set of first category. In  \S 7 we will show
that the concrete function $h$ given by (5.3) is not excluded (so, 
$F''$ is discontinuous).

      Let us consider the weighted shift operator $L_\beta$ in $A(\Omega)$ :

\[  (L_\beta U)(t)=\beta(t) U(ft). \]
Then we can rewrite (5.1) in the following way:
\begin{equation}
    (L_\beta -I)U=-h.
\end{equation}
This leads us to the problem of spectral properties of the weighted shift
operator.

  It is known  [27]  that the spectral radius 
$r_\beta$ of $L_\beta$ in the disk-algebra 
can be calculated by the formula:

\begin{equation}
        \ln r_\beta =\sup_{\nu\in M(f)}\chi_\nu(\beta),
\end{equation}
and the spectrum of $L_\beta$ is the unit disk:
${\rm spec}(L_\beta)=\{\lambda : \mid\lambda\mid\leq r_\beta     \}$.

If $U$ is an eigenfunction of $L_\beta$ then the function $\ln\mid\beta\mid$
is homologous to a constant:
\[\ln\mid\beta\mid=\ln\mid U(ft)\mid -\ln\mid U(t)\mid +\ln\mid\lambda\mid\]
which is not the case by Assumption (ii) above. So, the operator
$L_\beta -\lambda I$ is injective for all $\lambda$.By the Banach theorem on
 the inverse operator (see [28]), for
 $\lambda\in {\rm spec}(L_\beta)$ the image ${\rm Im}(L_\beta -\lambda I)$
is the set of first Baire category. Thus, the equation 
$L_\beta U-\lambda U=h$ for $\mid\lambda\mid\leq r_\beta$ and generic
$h\in A(\Omega)$ has no solutions in $A(\Omega)$.

    By (5.2), $r_\beta >1$, and hence $1\in{\rm spec L_\beta}$. So, the equation
(5.6) is non-solvable in $A(\Omega)$ for generic $h\in A(\Omega)$. For such
an $h$, the analytical solution (5.4) is not continuous up to the boundary,
as we have asserted.

    In conclusion  let us mention a rough obstruction for (5.1) to be
solvable related to non-invertibility of $f$. Let
$\overline t =(t_0,t_{-1},...)\in\overline\Omega$ be an inverse orbit of
$f$. Iterating (5.1) we get
\begin{equation}
 U(t)=\beta_{-n}(\overline t)U(t_{-n})-
    \sum_{k=1}^n{\beta_{-k}(\overline t)h(t_{-k})},
\end{equation}
where $\beta_{-k}(\overline t)=[\beta(t_{-1})...\beta(t_{-k})]^{-1}$.
 By the ergodic theorem

\begin{equation}
\lim_{n\rightarrow\infty}\frac{1}{n}\ln\beta_{-n}(\overline t)=
-\lim_{n\rightarrow\infty}\frac{1}{n}\sum_{k=1}^{n}\ln\beta(\overline f^{-k}t)=
        -\int{\ln\mid\beta\mid d\mu}<0
\end{equation}
for $\overline\mu$-almost all $\overline t\in\overline\Omega$.
For such a  $\overline t$ the sequence $\beta_{-n}(\overline t)$ exponentially
converges to zero, and we can consider the following 
$\overline\mu$-measurable function on $\overline\Omega$

\begin{equation}
   G(\overline t)=-\sum_{k=1}^\infty \beta_{-k}(\overline t)h(t_{-k})
\end{equation}
It follows from (5.8) that for continuous $h$ we have 
$G(\overline t)=U(t)$, i.e.\ the a priori multi-valued function $G$ turns  
out to be single-valued. It is the necessary condition for solvability of
the equation (5.1).

    Remark that this condition is almost sufficient. Namely, if the function
$G(\overline t)$ is single-valued, then it gives the holomorphic solution
in $\Omega^* =\Omega\verb+\+\{0\}$ which is continuous up to the boundary.
 However, this solution has a singularity at zero since $\beta(0)=0$.

\section{Asymptotics of $\cal F''$ along almost all geodesics 
 (conditional result)}

      In this section we assume that the function $U(t)$
given by (5.5) is not continuous up to the boundary. Under this assumption
we will calculate its asymptotics along almost all (in the sense of harmonic
measure) geodesics $B_\theta$. It turns out to be the
following:

\begin{equation}
 \overline{\lim_{r\rightarrow 1}}\,
 \frac{\ln\mid U(re^{i\theta})\mid }{-\ln(1-r)}=
 \frac{\chi_\mu(\beta)}{\chi_\mu}.
\end{equation} 

{\bf Remark.} By Makarov's theorem [31], 
$\ln l(t)\sim\ln(1-r)$ for $\mu$-almost all geodesics . So, (6.1) coincides
with the required asymptotics (1.5) (take in account formulas (4.4) and (5.4)
for characteristic exponents) 

     First let us give heuristic argument yielding (6.1).
 By (5.8),(5.10) we have

\begin{equation}
   \beta_{-n}(\overline t)U(t_{-n})\rightarrow U(t)-G(\overline t)
\end{equation}
for almost all $\overline t=(t,t_{-1},...)\in\overline\Omega $.
By (5.9),

     \[ \mid \beta_{-n}(\overline t)\mid \sim (\frac{2}{b})^n ,\]
where the sign $"\alpha_n\sim\beta_n"$ means

  \[ \lim_{n\rightarrow\infty}\ln\frac{\alpha_n}{\beta_n}=0.\]     
Hence,

\begin{equation}
   U(t_{-n})\sim(\frac{b}{2})^n.
\end{equation}

      Setting $t_{-n}=r_{-n}e^{i\theta_{-n}}, t=re^{i\theta}$, we get
$r_{-n}=r^{b^{-n}}$ ,and hence  

\[  b^n=\ln\frac{r}{r_{-n}}\sim\frac{1}{1-r_{-n}}\]
(since $r$ is fixed). Comparing this with (6.3), we find

  \[ U(t_{-n})\sim(\frac{1}{1-r_{-n}})^{\alpha_c} \]
where $\alpha_c=\chi_\mu(\beta)/\chi_\mu$ is the critical exponent written above.

    The main shortcoming of this calculation is that the points $t_{-n}$
don't lie near a single radius $B_\theta$ but are wandering along
$\partial\Omega$. However, they quite often approach almost every geodesic
which allows to turn our argument into rigorous.

    Let us consider the multi-valued function $G(\overline t)$ on the 
covering space $\overline\Omega$ given by the series (5.10). This function
is correctly defined due to (5.2), and satisfies the following equation

\begin{equation}
  \beta(t) V(\overline f\overline t)-V(\overline t)=-h(t),
\end{equation}
the lift of equation (5.1) on $\overline\Omega$. In order to substantiate
the asymptotics we need some analytical properties of the $G$.

  Let us put on $\overline\Omega\simeq [0,1]\times\Sigma_b$  measure $\lambda$,
the product of the Lebesgue measure on [0,1] and the Bernoulli measure
$\overline\mu$ on $\Sigma_b$. Consider the space $A(\overline\Omega)$ of functions
measurable on $\overline\Omega$ , analytic on almost all Riemann
surfaces $W(\overline\varepsilon_-)$ and continuous up to the boundary 
in their inner topology (see \S 4). 
Remark that they probably are not defined at 0.

\begin{Lemma}
The function $G(\overline t)$
given by {\rm (5.10)} belongs to the space $A(\overline\Omega)$
\end{Lemma}
{\bf Proof.} By the ergodic theorem,
 $\beta_{-n}(\mid\overline t\mid)\sim e^{-an}\quad (n\rightarrow\infty)$
for $\overline\mu$-almost all $\overline t=(1,\overline\varepsilon)\in\Gamma.$
Applying  (4.6) to the function $\rho=-\ln\mid\beta\mid$, we get the
similar asymptotics for all 
$\overline\tau\in W(\overline t)\equiv W(\overline\varepsilon_-)$.
Moreover, for all $\theta >0$ there exists $Q(\delta,\theta)$ such that

\[\mid\beta_{-n}(\overline\tau)\mid\leq
e^{C(\delta)}\mid\beta_{-n}(\overline t)\mid\leq
Q(\delta,\theta)e^{-(a-\theta)n},\]
provided $\mid\overline\tau\mid \geq\delta >0$. Consequently,the series (5.10)
converges uniformly on $W_\delta(\overline t)$ which yields the required    
statement. \QED

  There is the Bernoulli measure $\overline\mu$ on the solenoids 
$\overline\Gamma_r\simeq\Sigma_b$, and one can consider the corresponding
spaces $L^\kappa(\overline\Gamma_r),\quad\kappa >0.$
Now let us prove the Main technical lemma.

\begin{Lemma}
  For sufficiently small $\kappa >0$ and any $r\in(0,1]$, the function
$G(\overline t)$ belongs to the space $L^\kappa(\overline\Gamma_r)$:
\[\int_{\Sigma_b}{\mid G(r,\overline\varepsilon) \mid^\kappa}
d\overline\mu(\overline\varepsilon) <\infty.\]
Moreover, for any $\delta >0$ the same is true for the function
\[G_\delta(\overline\varepsilon)=\sup_{r:\delta\leq r\leq 1}
\mid G(r,\overline\varepsilon)\mid . \]   
\end{Lemma}
{\bf Proof.} Applying  (4.6) to the function $\rho=\ln\mid\beta\mid$,
we get for $r\geq\delta$:
\[\mid\beta_{-n}(r,\overline\varepsilon)\mid\leq
Q(\delta)\mid\beta_{-n}(1,\overline\varepsilon)\mid.\]

Consequently,

\[\mid G_\delta(\overline\varepsilon)\mid\leq
(\mid\mid\psi\mid\mid_\infty Q(\delta))^\kappa
\sum_{n=1}^{\infty} \mid\beta_{-n}(1,\overline\varepsilon)\mid^\kappa ,  \]

and further

\[\int_{\overline\Gamma}{\mid G_\delta\mid^\kappa d\overline\mu}\leq
(\mid\mid\psi\mid\mid_\infty Q(\delta))^\kappa
\sum_{n=1}^{\infty} \int_{\overline\Gamma}{\mid\beta_{-n}\mid ^\kappa 
d\overline\mu.} \]

But
\[\int_{\overline\Gamma}{\mid\beta_{-n}(\overline t)\mid ^\kappa 
d\overline\mu(\overline t)}
=\int_{\overline\Gamma}{\mid\beta_{n}(\overline f^{-n}\overline t)\mid ^{-\kappa}
 d\overline\mu(\overline t)}
=\int_{\overline\Gamma}{\mid\beta_{n}(\overline t)\mid ^{-\kappa} 
d\overline\mu(\overline t)}
=\int_{\Gamma}{\mid\beta_n (t)\mid ^{-\kappa} d\mu(t)} \]

Thus, it is enough to show that for sufficiently small $\kappa >0$ the last
integral exponentially tends to zero as $n\rightarrow\infty$. To this end
we will apply thermodynamical formalism.

   Let $\nu_\kappa$ be the Gibbs measure corresponding to the potential
\mbox{$-\kappa\ln\mid\beta\mid$},and $P_\kappa =P_f (-\kappa\ln\mid\beta\mid) $
be the corresponding pressure. By (4.1), 

\[\mid\beta_n(t_{\varepsilon_0 ...\varepsilon_{n-1}})\mid^{-\kappa}
\asymp\nu_\kappa[\varepsilon_0 ...\varepsilon_{n-1}]\exp(nP_\kappa)\]
for any point $t_{\varepsilon_0 ...\varepsilon_{n-1}}$ from the cylinder
$[\varepsilon_0 ...\varepsilon_{n-1}]$. 

  As $\mu[\varepsilon_0 ...\varepsilon_{n-1}]=b^{-n}$ then

\[\int{\mid\beta_n\mid ^{-\kappa} d\mu}\asymp
 b^{-n}\sum_{1\leq\varepsilon_i\leq n}
{\nu_\kappa[\varepsilon_0 ...\varepsilon_{n-1}]}\exp(nP_\kappa)=
\exp[n(P_\kappa-\ln b)].
\]

   Hence, it is sufficient to check that for small enough $\kappa >0$

\begin{equation}
   P_\kappa -\ln b <0
\end{equation}

But $P_0 =h(f) =\ln b$, and by (4.3) and (5.2) we have

\[ \frac{dP_\kappa}{d\kappa}\mid_{\kappa =0} =
  -\int{\ln\mid\beta\mid d\mu}<0.
\]

\psfig{figure=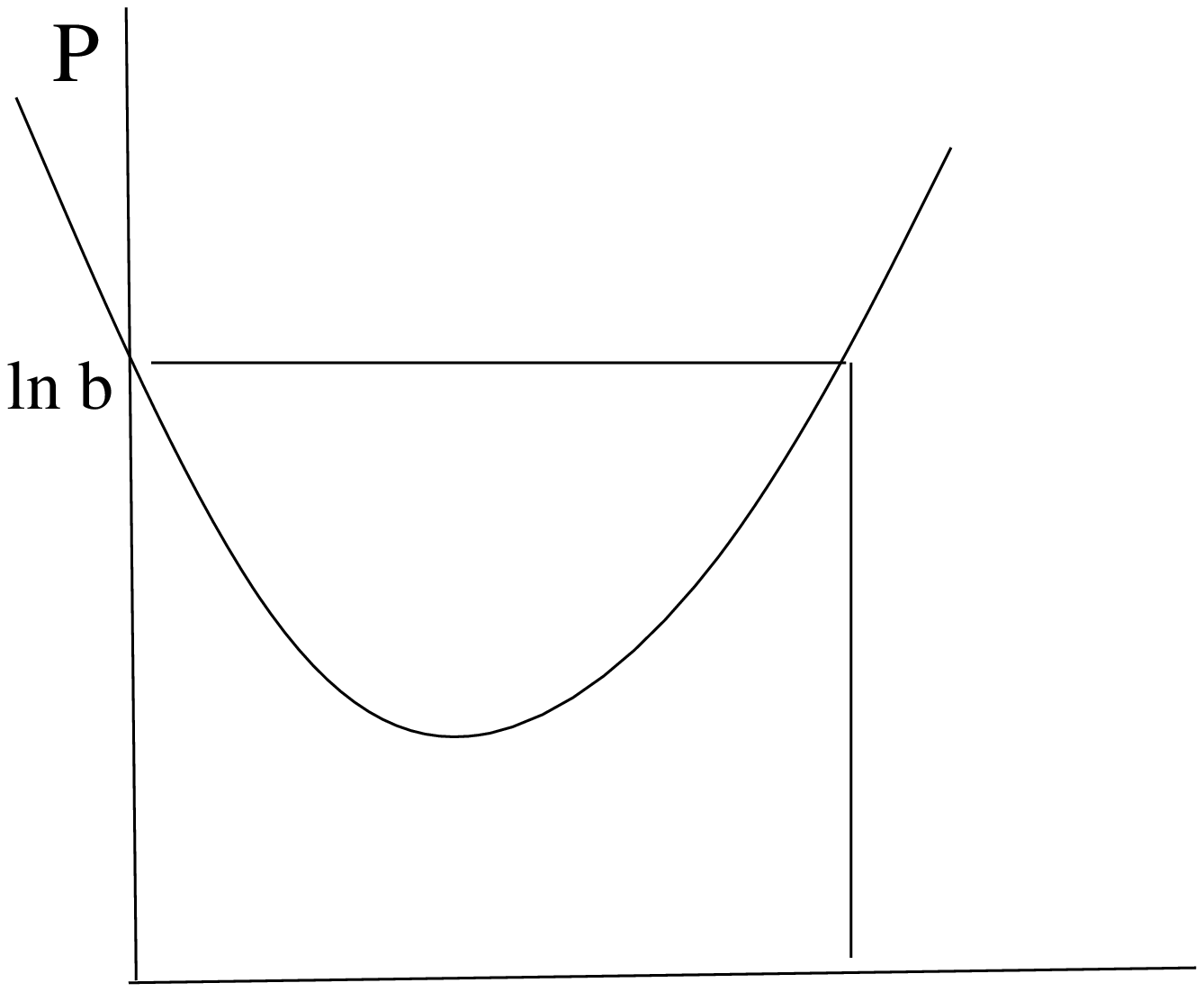,height=1.8in,width=3in}
\vspace{.4in}
\hspace{2in} Figure 4
\vspace{.4in}

Hence, (6.5) holds for any $\kappa\in(0,\kappa_0)$ (Fig. 4), and the lemma is
proved.   \QED

      So, we have two solutions of the equation (5.1), 
$U(\overline t)\equiv U(t)$ given by (5.4) and $G(\overline t)$ given
by (5.10). The first of them comes from the single-valued function
on $\Omega_0$; the second is multi-valued on $\Omega_0$ but possesses
better boundary properties. The difference $V(\overline t)=U(t)-G(\overline t)$
satisfies the homogeneous equation

\begin{equation}
    \beta (t) V(\overline f\overline t)= V(\overline t).
\end{equation}
Since $U(t)$ is not continuous up to the boundary (see the beginning of
 the section), while $G\in A(\overline\Omega)$, we have $V\not =0$.
We are going to find the asymptotics of $V$ along $\overline\mu -$almost all
"geodesics" $\overline B_{\overline\varepsilon}=\{(r,\overline\varepsilon)\in
\overline\Omega : 0\leq r\leq 1\}$ as $r\rightarrow 1$. Since $G$ is
continuous on almost all sheets, the same asymptotics will be valid for $U$.

   Set
\[V_0(\overline\varepsilon)=\max_{1/2\leq r\leq b/2} 
\mid V(r,\overline\varepsilon)\mid.\]
By Lemma 6.2, $V_0\in L^\kappa$. Hence

\[\int_{\Sigma_b}{\ln \mid V_0(\overline\varepsilon)\mid
 d\overline\mu(\overline\varepsilon)} < +\infty.\]
Now we can apply the ergodic theorem. It follows

\begin{equation}
\overline{\lim_{n\rightarrow\infty}}\quad\frac{1}{n}\ln\mid V_0
(\sigma^n\overline\varepsilon)\mid\leq 0
\end{equation}
for almost all $\overline\varepsilon\in\Sigma_b$.
Let us consider a geodesic $\overline B_{\overline\varepsilon}$
 corresponding to such an $\overline\varepsilon$.  
Let $r\geq 1/2$ and $n=n(r)$ is such a number that
\begin{equation}
\frac{1}{2}\leq r^{b^n}<(\frac{1}{2})^{1/b}
\end{equation}

  It follows from (6.6) and (6.7) that

\[ \mid V(r,\overline\varepsilon)\mid= \mid\beta_n(r,\overline\varepsilon_+)\mid
\mid V(r^{b^n},\sigma^n\overline\varepsilon)\mid
\prec \mid\beta_n(r,\overline\varepsilon_+)\mid
\]
where the sign $"\alpha\prec\beta"$ means 
 $\overline{\lim}\frac{1}{n}\ln (\beta_n/\alpha_n) \leq 0$. 
Since by the ergodic theorem,
$\mid\beta_n (r,\overline\varepsilon_+)\mid\sim e^{na}$
for almost all $\overline\varepsilon_+$, we conclude

\[\overline{\lim_{r\rightarrow 1}}\, \frac{1}{n(r)} 
\ln V(r,\overline\varepsilon)
\leq \chi_\mu(\beta).
\]

From (6.8) we find

\[n(r)=-\frac{\ln(1-r)}{\ln b} +O(1).\]

Two last estimates yield
\begin{equation}
\overline{\lim_{r\rightarrow 1}} \,
\frac{\ln V(r,\overline\varepsilon)}{-\ln(1-r)}\leq\frac{\chi_\mu(\beta)}{\ln b}.
\end{equation}

In order to get the opposite estimate let us consider the set
\[X_\delta =\{\overline\varepsilon\in\Sigma_b :
V(1/2,\overline\varepsilon)\geq\delta\}\].
 As $V\not\equiv 0$, $\lambda(X_\delta)>0$
for sufficiently small $\delta>0$. Hence, almost all orbits 
$\{\sigma^n\overline\varepsilon\}_{n=0}^\infty$ pass through $X_\delta$
infinitely many times. For such an $\overline\varepsilon$ let us consider
a sequence $n(k)\rightarrow\infty$ for which 
$\sigma^{n(k)}\overline\varepsilon\in X_\delta$. 
Set $r_n=(1/2)^{b^{-n}}$. Then (6.6) implies

\[ \mid V(r_{n(k)},\overline\varepsilon)\mid=
\mid \beta_{n(k)}(r_{n(k)},\overline\varepsilon_+)\mid
\cdot\mid V(1/2,\sigma^{n(k)}\overline\varepsilon)\mid\geq
\delta\mid \beta_{n(k)}(r_{n(k)},\overline\varepsilon_+)\mid.
\]
This implies the inequality opposite to (6.9). Thus, for almost all
$\overline\varepsilon\in\Sigma_b$ we have
\[ \overline{\lim_{r\rightarrow 1}}\quad
\frac{\ln V(r,\overline\varepsilon)}{-\ln(1-r)}
=  \frac{\chi_\mu(\beta)}{\ln b},\]
and the required asymptotics (6.1) is proved.

\section{$\cal F''$ is not continuous up to the boundary}

Let us prove first that one of the derivatives is discontinuous  
on the interval $[0,t_c]$ at the real
critical point $t_c$ (it is a critical point from the thermodynamical viewpoint;
from the dynamical viewpoint, $t_c$ is a repelling fixed point ).
Denote by $\lambda=f^n(t_c) >1$ the multiplier of this
point.

\begin{Lemma}
 Let $l$ be a natural number for     which $\lambda^l > 2b$.
Then $f^{(l)}$ is discontinuous on the interval $[0,t_c]$ at $t_c$.
\end{Lemma}
{\bf Proof.} Let us linearize $f$ at $t_c$:
\[\psi (fz)=\lambda\psi (z)\]
where $\lambda$ is an analytic function in an neighbourhood of $t_c$,
$\psi (t_c)=0, \psi '(t_c)=1$ (the {\em K\"{o}nigs function}).
Clearly, $\psi$ can be analytically continued on the interval $(0,t_c]$, and 
 one-to-one maps it onto the axis $(-\infty, 0]$.

 Set

\[\stackrel{\sim}{F} =F\circ\psi^{-1},\hspace{.5in}
 \stackrel{\sim}{g} =g\circ\psi^{-1}.\] 

Then $\stackrel{\sim}{F}$ satisfies the following functional equation:

\[\frac{1}{2b}\stackrel{\sim}{F}(\lambda z) -\stackrel{\sim}{F}(z)
=\stackrel{\sim}{g}(z).\]

Due to the linearization we immediately get the functional equation for 
$\stackrel{\sim}{F}^{(l)}$:

\[\frac{\lambda^{(l)}}{2b}\stackrel{\sim}{F}^{(l)}(\lambda z) -\stackrel{\sim}{F}^{(l)}(z)
=\stackrel{\sim}{g}^{(l)}(z).\]

If $\stackrel{\sim}{F}^{(l)}$ is continuous on the semi-axis $(-\infty,0]$,
then it can be given on it by the following series (compare (5.10)):
\[\stackrel{\sim}{F}^{(l)}=-\sum_{k=1}^{\infty}
\frac{(2b)^k}{\lambda^{nk}}\stackrel{\sim}{g}^{(l)}(\frac{z}{\lambda^k})\]
which gives the analytical continuation of $\stackrel{\sim}{F}^{(l)}$
through 0. Hence $F$ can be analytically  continued through $t_c$
contradicting Lemma 3.2.    \QED

{\bf Remark.} The same argument can be applied to any fixed point 
$\alpha\in\partial\Omega$: some derivative of $F\mid\Omega^\circ$
 must be discontinuous at $\alpha$.

\begin{Lemma}
 The function $F''$ is not continuous up to the boundary of $\Omega$.
\end{Lemma}
{\bf Proof.}Assuming the reverse, we will show
that all derivatives $F^{(n)}$ should be continuous up to the boundary,
contradicting Lemma 7.1. 

Let us consider the functional equations for the derivatives of $F$
(cf. (3.6) and (5.3)):

\begin{equation}
\frac{(f')^{n+1}}{2b}F^{(n+1)}\circ f-F^{(n+1)}=-h_{n}
\end{equation} 
where $h_{n}$ can be expressed via the derivatives of $F$ of order $\leq n$:

\begin{equation}
h_n=-g^{(n)}-\frac{1}{2b}\sum_{k=1}^{n-1}
(\frac{d(f')^k}{dt}\cdot\frac{d^k F}{dt^k}\circ f)^{(n-k-1)}
\end{equation}
It is convenient to use the following metric in $\Omega^0$:
\[d(z,\zeta)=\inf l(\gamma)\]
where inf is taken over all rectifiable paths $\gamma$ connecting $z$ and
$\zeta$.

Assume by induction that all derivatives $F^{(k)}, k=1,2,...,n$, 
are continuous in $\Omega$ (n=2 is the base of induction). 
Then the derivatives $F^{(k)}, k=1,2,...,n-1$, should be Lipschitz continuous
with respect to the metric d. By (7.2), $h_n$ should possess the same property. 
But if $F^{(n)}$ is continuous in $\Omega$, $n\geq 2$, then it can be given
by the series (5.10):
\begin{equation}
F^{(n)}(t)=-\sum_{k=1}^{\infty}\frac{(2b)^k}
{[(f^{\circ k})'(t_{-k})]^n}h_{n-1}(t_{-k})
\end{equation}
for a $\mu$-typical inverse orbit $\overline t=\{t_0,t_{-1},...\}$.
It follows from here (taking in account the expanding property of $f\mid\Gamma $)
that $F^{(n)}$ is also Lipschitz continuous with respect to d. Hence, 
$F^{(n+1)}$ is bounded in $\Omega^0$. It is enough for $F^{(n+1)}$ to be given
by the series (5.10) and, hence, to be continuous up to the boundary.\QED

\section{Concluding remarks}
Note that the complex critical exponent calculated in this paper differs
from the usual real critical exponent $\alpha$  at $t_c$. 
In fact, there is the general scheme including both cases.
Namely, one can
associate to any invariant measure $\nu$ on $\partial\Omega_0$ its own
critical exponent $\alpha_\nu$, i.e.\ the exponent of power growth of an
appropriate derivative
of the free energy along $\nu$-typical geodesics. 
More specifically, let
\[m_\nu=[\frac{\ln 2b}{\chi_\nu}]+1,\; \alpha_\nu=1-\{\frac{\ln 2b}{\chi_\nu}\} \]
where  $[a]$ and $\{a\}$ denote the entire and the fractional part of $a$
respectively, and $\chi_\nu$ is the characteristic exponent of $\nu$. Then  
the following general formula should be true:

\begin{equation}
  \lim_{t\rightarrow\tau ,\, t \in{B_\tau}}  \frac{\ln \mid
{\cal F}^{(m_\nu)}(t)\mid}{-\ln l(t)}=\alpha_\nu
\end{equation}
for $\nu$-almost all geodesics.

  Indeed, following the scheme of the present paper, we should find
the first
$m$ for which the function $\beta=(f')^m/2b$ has positive characteristic
exponent:
\[\chi_\nu(\beta)=m 
\chi_\nu-\ln 2b>0,\]
i.e.\ $m>\ln 2b/\chi_\nu$ (assume that $\ln 2b/\chi_\nu$ is non-integer).
Then by formula (6.1) we get (8.1).

For the $\delta$-measure  concentrated at the critical point $t_c$ we
obtain the usual formula of the renorm-group theory (see [29] , [13]): 
\[\alpha=1-\{\frac{\ln 2b}{\ln f'(t_c)}\}=
1-\{\frac{d\ln 2}{\ln f'(t_c)}\}.\]
The similar formula holds for any periodic point on $\partial\Omega$
(without any changes in the proof).
 
   The main technical problem in proving (8.1) for general $\nu$
 is related to the fact
that if $m>2$ then the right-hand side of the equation (7.1) for $F^{(m)}$
is discontinuous on the boundary. The same problem arises if one wants
to calculate the critical exponent for the measure of maximal entropy on
the boundary of of the "high-temperature" basin $\Omega_1$ .
The true formula should be $\alpha=1-\{\ln b/\ln 2\}.$
It is interesting also to find the complex critical exponent for $b=d=2$.

    The other problem is to study the global properties of free energy on the
whole Riemann sphere.  
 They have to do with Gibbs measures on $J(f)$.

We finish the paper with the following important remark.
 One can show that the free energy $F$ can be represented as the logarithmic
potential of the measure of maximal entropy  of $f$. This gives another 
approach to the circle of problems under consideration and a nice relation
of the critical exponent to the local dimension of the measure of maximal
entropy. We are grateful to P. Moussa and A. Eremenko for interesting 
discussions of this point.

\end{document}